\documentclass[11pt,a5paper,twoside]{article}
\usepackage[a5paper,margin=13mm]{geometry}
 
 \title{A multiscale scheme accurately simulates macroscale shocks in an equation-free framework}
 \author{John~Maclean\thanks{\protect\raggedright%
School of Mathematical Sciences, University of Adelaide, South Australia.
\url{http://www.adelaide.edu.au/directory/john.maclean}}
\and J.~E.~Bunder\thanks{\protect\raggedright%
School of Mathematical Sciences, University of Adelaide, South Australia.
\protect\url{mailto:judith.bunder@adelaide.edu.au},
\protect\url{http://orcid.org/0000-0001-5355-2288}}
\and A.~J.~Roberts\thanks{\protect\raggedright%
School of Mathematical Sciences, University of Adelaide, South Australia.
\url{http://www.maths.adelaide.edu.au/anthony.roberts},
\url{http://orcid.org/0000-0001-8930-1552}}
\and 
I.~G. Kevrekidis
\thanks{Departments of Chemical and Biomolecular Engineering and Applied Mathematics and Statistics, Johns Hopkins University, Baltimore, Maryland, USA.
\protect\url{https://orcid.org/0000-0003-2220-3522}}
}

\date{\today}

\usepackage[T1]{fontenc}
\usepackage[dvipsnames]{xcolor}
\usepackage{enumitem}

\usepackage{pgfplots,ifthen} 
\pgfplotsset{compat=newest} 
\usepgfplotslibrary{external} 
\usetikzlibrary{decorations.markings}
\usetikzlibrary{shapes,arrows,fit}
\usetikzlibrary{positioning}
\usetikzlibrary{calc} 

\usepackage{natbib}
\AtEndDocument{%
    
    \bibliography{ajr,bib,bibexport}
    }

\usepackage{url,microtype,amsmath,amssymb,graphicx,defns,hyperref,doi}
\hypersetup{colorlinks,allcolors=RoyalBlue,breaklinks,%
            bookmarks,bookmarksopen}
\makeatletter
\AtBeginDocument{{\def\and{, }\def\thanks#1{}%
  \hypersetup{
    pdfauthor={\@author},
    pdftitle={\@title}}}
    }
\makeatother

\usepackage{mycleveref}

\Vec f\Vec x\Vec u \Vec X \Vec U
\Bb R  \Bb X

\graphicspath{{figs/}}

 \definecolor{Mcol}{RGB}{217,  95,   2}
\definecolor{Mtxt}{RGB}{165,  71,   2} 
 \definecolor{mcol}{RGB}{0,  123,   123} 
\definecolor{mtxt}{RGB}{0,  100,   100} 
 \definecolor{excol}{RGB}{76,  0,   153}  
\definecolor{extxt}{RGB}{80,  0,   163}  

\begin{document}

\maketitle

\begin{abstract}
Scientists and engineers often create accurate, trustworthy, computational simulation schemes---but all too often these are too computationally expensive to execute over the time or spatial domain of interest.
The equation-free approach is to marry such trusted simulations to a framework for numerical macroscale reduction---the patch dynamics scheme.
This article extends the patch scheme to scenarios in which the trusted simulation resolves abrupt state changes on the microscale that appear as shocks on the macroscale.
Accurate simulation for problems in these scenarios requires extending the patch scheme by capturing the shock within a novel patch, and also modifying the patch coupling rules in the vicinity in order to maintain accuracy.
With these two extensions to the patch scheme, straightforward arguments derive consistency conditions that match the usual order of accuracy for patch schemes.
The new scheme is successfully tested on four archetypal problems.
This technique will empower scientists and engineers to accurately and efficiently simulate, over large spatial domains, multiscale multiphysics systems that have rapid transition layers on the microscale.
\end{abstract}

\section{Introduction}

The modeling of scientific and engineering phenomena is often complicated by the presence of fast, fine-scale processes entangled with the long-lasting, macroscale, system-wide, processes that are of interest.
This article contributes to developing the so-called ``equation-free'' methodology for efficient system level simulation of such complex multiscale phenomena \cite[e.g.]{Kevrekidis03b, Kevrekidis09a, Sieber2018}.
This methodology accurately predicts the macroscale, system-wide, coarse scales but only requires the fast, fine-scales to be resolved  on small patches and/or in small bursts of the space-time domain.
Projective integration uses short bursts of the full microscale simulation in time to learn the macroscale information that empowers long-time prediction \cite[e.g.]{Gear02b, Erban2006, Givon06}.
To complement projective integration, we focus on problems that have multiple scales in space by further developing the patch scheme which computes on only small microscale patches in space and yet makes accurate macroscale prediction \cite[e.g.]{Gear03, Hyman2005, Samaey03b, Samaey04, Roberts06d}.
Specifically, here we address and resolve issues that arise when the effective macroscale dynamics exhibits localised `shocks', which nonetheless the microscale simulation could resolve except for the resultant computational expense.  

The detailed structure of shocks, cracks, grain boundaries, and dislocations, are only resolved on a micro-scale (such as an atomic simulation), so any description of shock\slash crack development over large-scales is a multiscale issue.  
The importance and dynamic nature of such shock propagation \cite[e.g.]{Takayama04} and the resulting discontinuities in material structure have led to many mathematical approaches \cite[e.g.]{Guenter2011, Luskin2013}, including the quasi-continuum methodology of \cite{Tadmor1996} which has some features cognate to our approach.
In its hybrid design, the patch scheme (\cref{sec:pat}) that we develop further here is ideally suited to multiscale modelling of systems with complex microscale detail \cite[]{Kevrekidis09a, Samaey08, Roberts2014a}. 

\begin{figure}
    \centering
\begin{tabular}{cc}
poor attempt
&accurate simulation
\\
    \includegraphics[]{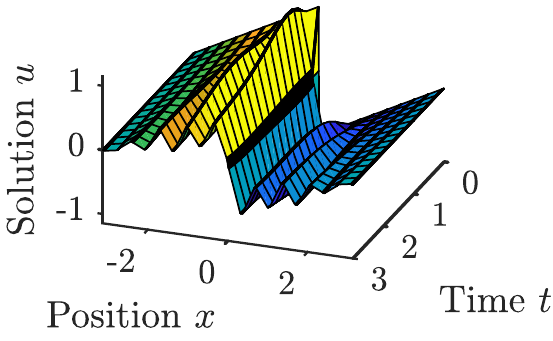}
&   \includegraphics[]{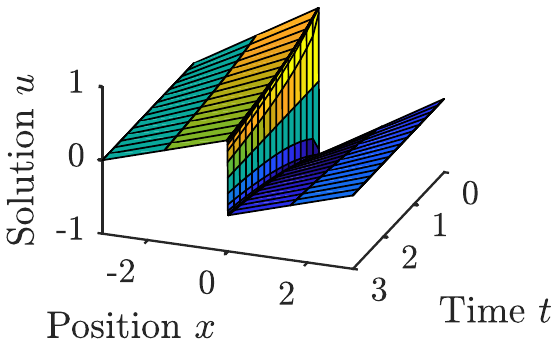}
\end{tabular}
    \caption{
Two attempts to capture a shock by the patch scheme. 
The black band in the centre of both graphs is a micro-scale discretisation consisting of \(25\)~mesh points uniformly spaced in a `patch' of width~\(0.05\), and this patch is appropriately placed to resolve the shock. 
Nonetheless, the left-hand straightforward attempt fails to capture the shock, whereas the right-hand graph shows that \cref{sec:dpat}'s modifications makes excellent predictions.
}
    \label{fig:stdPat}
\end{figure}

The key idea of the patch scheme is to replace an expensive simulation across a `large' domain with a set of simulations in much smaller, well-separated, patches in the domain (\cref{sec:pat}). 
With appropriate coupling between the patches macroscale predictions are provably accurate \cite[e.g.]{Roberts06d, Bunder2013b}, and so the scheme implicitly performs a macroscale reduction.
Furthermore, in this equation-free approach, such an effective macroscale reduction is obtained blindly (almost); no analysis of the micro-scale structures is necessary, only knowledge that the microscale model is accurate. 
The scheme is cognate to \emph{computational homogenization} \cite[e.g.]{Geers2010, Saeb2016, Geers2017} and to \emph{numerical homogenization} \cite[e.g.]{Craster2015, Owhadi2015, Peterseim2019, Maier2019}. 
In diverse problems, the patch scheme has been shown to be accurate to order~$\Ord{H^{2\Gamma}}$, where $H$~is the distance between patches and~$\Gamma$ is the order of coupling between patches  \cite[e.g.]{Roberts06d, Roberts2011a, Cao2014a}, including when the microscale is heterogeneous \cite[e.g.]{Bunder2013b}. 
The documentation of the \emph{Equation-Free Toolbox}\footnote{\url{https://github.com/uoa1184615/EquationFreeGit} for \script.} gives a practical, contemporary introduction to the patch scheme and includes diverse examples.  
However, the patch scheme as developed previously does \emph{not} capture shocks, as illustrated by the poor simulation shown in \cref{fig:stdPat}(left).

This article begins to develop a new extension to the patch scheme in order to accurately simulate systems including one or more shocks (\cref{sec:dpat}).
In its shock-resolving capability, the extension has some advantages in comparison to other approaches that aim to resolve a system with a shock.
One classical approach is the Rankine--Hugoniot (\textsc{rh}) conditions that are used to relate solutions on either side of the shock---but in some application areas these conditions do not apply.
We demonstrate that our patch scheme accurately solves problems (\cref{exUstat,exUform}) for which the \textsc{rh} conditions do not apply.
Sometimes numerical formulations are employed to simulate systems with shocks, such as adaptive mesh methods \cite[e.g.]{Huang10} or sparse grids \cite[e.g.]{JakemanArchibaldXiu11}.
Neither of these two methods perform a reduction from micro- to macro-scale along the lines of the patch scheme, so they are best for a different class of problems.

The numerical examples in this article (\cref{sec:modp}) are all archetypes for the extended patch scheme, and should be understood in the following context. 
We consider the key complication to address is a single shock in a simulation that is otherwise appropriately simulated by the patch scheme.
Therefore our canonical problems do not contain complicated dynamics except at the shock---it is sufficient to establish that the new extension can accurately resolve a shock.

\section{Four archetype problems}
\label{sec:modp}

In this first development of shocks in the patch scheme we restrict attention to \pde\ systems in 1D space, but the techniques would also apply immediately to multi-D space, of large extent in one direction, when the spatial patches extend across the cross-sections \cite[e.g.]{Alotaibi2017a}.
As archetype problems, this article specifically seeks to predict the dynamics of a field~\(u(x,t)\) satisfying a modified form of Burgers' \pde, in 1D space, with diffusion coefficient~$\eps(u)$, namely
\begin{align} \label{be}
\D tu + u \D xu &= \eps(u) \DD xu \,.
\end{align}
We solve this \pde\  on the non-dimensional spatial domain \(-\pi<x<\pi\)\,, with an initial condition that $u(x,0)=u_0(x)$, and with Dirichlet boundary conditions of $u=u_L(t),u_R(t)$ at $x=\pm\pi$ respectively.

This \pde\ lets us consider a series of archetype problems, each exhibiting a shock, with small qualitative differences in the difficulty of simulation. 
When choosing a constant diffusivity $\eps(u)=\eps$ we recover Burgers' \pde, and we use this in problems~\cref{exStat,exForm} to generate simple shocks that demonstrate our novel approach in two well understood problems. 
One recent example application, among many, would be to the shock fronts in Burgers' \pde\ that model flow through complex networks affected by abrupt local changes \cite[e.g.]{Mones2014}.
Then, in problems~\cref{exUstat,exUform}, we choose a variable diffusivity~$\eps(u)$ that renders the \textsc{rh}~conditions inapplicable, so the properties of the shock cannot be analytically calculated. 

\begin{figure}
    \centering
\begin{tabular}{cc}
    \includegraphics[]{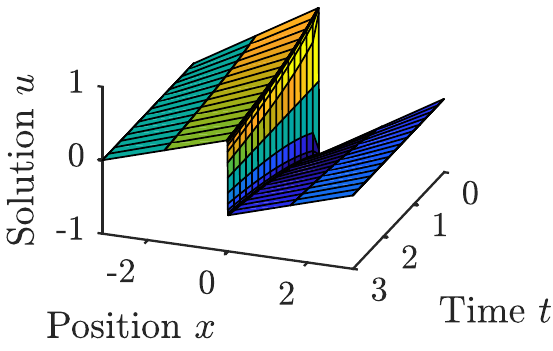}
&   \includegraphics[]{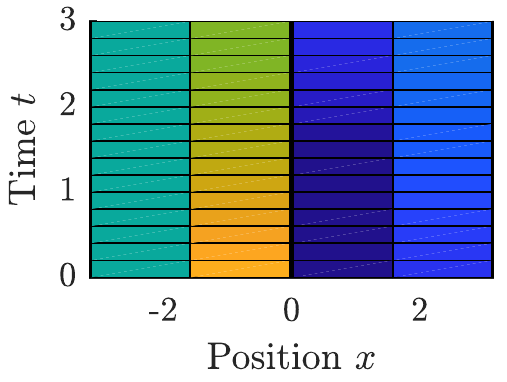}
\end{tabular}
    \caption{A numerical quadrature solution of~\cref{exStat}.
The left figure plots~$u(x,t)$ on the vertical axis, while the right figure shades each cell according to the values of~$u(x,t)$ at the cell vertices.
The solution is computed on \(25\)~tightly spaced points between~$-0.025$ and~$0.025$, and on four other points in the spatial domain.
The narrow band of points show the behaviour of the solution around the shock, while the coarse discretisation captures the (macroscale) behaviour away from the shock.  }
    \label{fig:tanh}
\end{figure}

\begin{figure}
    \centering
\begin{tabular}{cc}
    \includegraphics[]{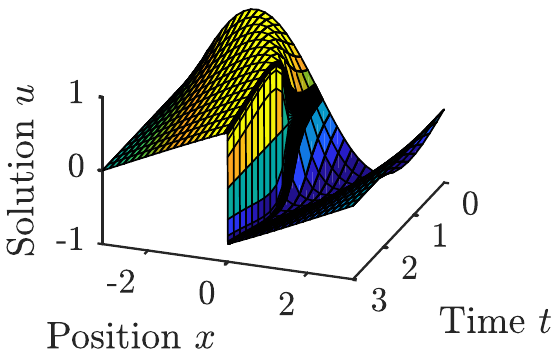}
&   \includegraphics[]{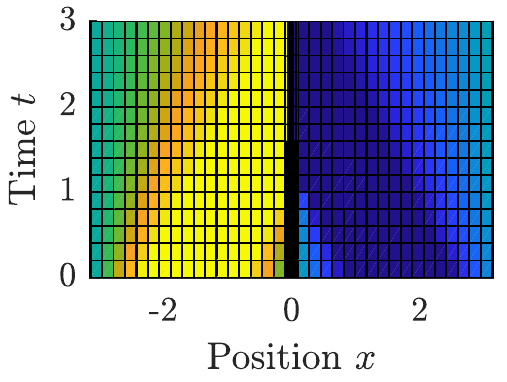}
\end{tabular}
    \caption{A numerical quadrature solution of~\cref{exForm}.
The left figure plots~$u(x,t)$ on the vertical axis, while the right figure shades each cell according to the values of~$u(x,t)$ at the cell vertices.
The solution is computed on \(100\)~tightly spaced points between~$-0.1$ and~$0.1$, and on \(34\)~other points in the full spatial domain.
The solution is initially smooth, but forms a shock at roughly $t = 1.2$. }
    \label{fig:sin}
\end{figure}

\begin{figure}
    \centering
\begin{tabular}{cc}
    \includegraphics[]{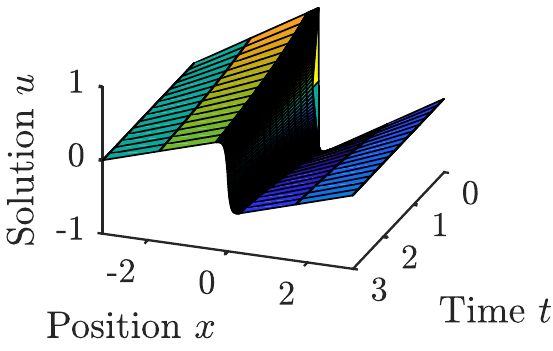}
&   \includegraphics[]{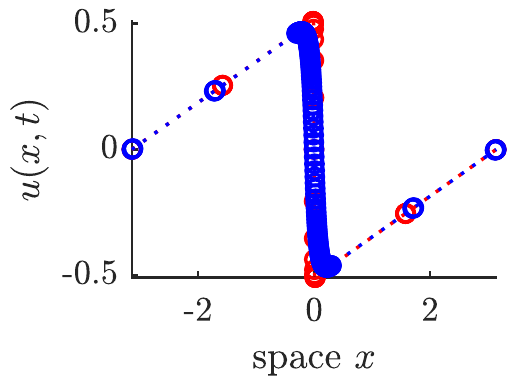}
\end{tabular}
    \caption{An accurate numerical simulation of~\cref{exUstat}.
The left graph plots~$u(x,t)$ on the vertical axis, while the right graph compares the simulation at time $t=3$ to the solution obtained for~\cref{exStat} at that time. 
The simulation is computed on \(1600\)~evenly spaced points in space, and the numerical time step is $5\cdot10^{-6}$. 
Then that fine scale solution is evaluated at spatial locations that are the patches in \cref{sec:dpat}.}
    \label{fig:Utanh}
\end{figure}

\begin{figure}
    \centering
\begin{tabular}{cc}
    \includegraphics[]{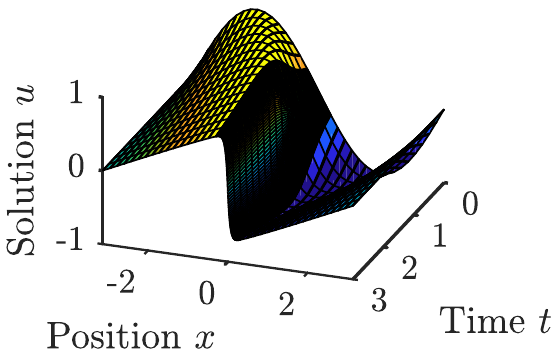}
&   \includegraphics[]{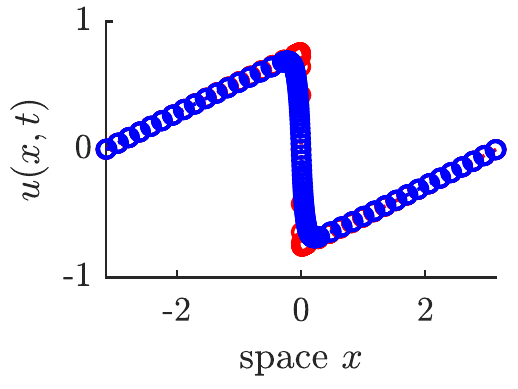}
\end{tabular}
    \caption{An accurate numerical simulation of~\cref{exUform}.
The left figure plots~$u(x,t)$ on the vertical axis, while the right figure compares the simulation at time $t=3$ to the solution obtained for~\cref{exStat} at that time.
The simulation is computed on \(1600\)~evenly spaced points in space, and the numerical time step is~$5\cdot10^{-6}$; then that fine scale solution is evaluated at spatial locations that are the patches in \cref{sec:dpat}.}
    \label{fig:Usin}
\end{figure}

We consider the following four archetype problems.
\begin{enumerate}[label=\textbf{M\arabic*}]
    \item \label{exStat} Fix the diffusion coefficient $\eps(u):= 0.001$.
Boundary conditions are that $u=0$ at $x=\pm\pi$.
The initial condition that
\begin{equation*}
    u_0(x) := \frac{{x}/{\pi} -\tanh\left({2x}/{\eps}\right)}{\tanh\left({2\pi}/{\eps}\right)}
\end{equation*}    
includes a rapid transition of width proportional to~$\eps$ centred at $x=0$\,, that appears as a shock on the macroscale.
\cref{fig:tanh} plots an exact solution.

    \item \label{exForm} Fix the diffusion coefficient $\eps(u) := 0.001$.
The boundary conditions are that $u=0$ at $x=\pm\pi$.
The initial condition is the smooth $u_0(x) := -\sin x$\,.
The solution, at about time $t=1.2$\,, forms a rapid transition at $x=0$ that appears as a macroscale shock (\cref{fig:sin}).

    \item \label{exUstat} Let the diffusion coefficient be the nonlinear $\eps(u) := 0.001 + 0.05|u|$.
Boundary and initial conditions as for~\cref{exStat}.
We chose this diffusion coefficient so that then \pde~\cref{be} cannot be written in conservation form, and hence the Rankine--Hugoniot conditions do not predict the shock speed nor relate the solutions on either side of the shock.
\cref{fig:Utanh} shows an accurate simulation, and contrasts it to the solution of~\cref{exStat}.

    \item \label{exUform} Set the diffusion coefficient to be the nonlinear $\eps(u) := 0.001 + 0.05|u|$.
Boundary and initial conditions are as for~\cref{exForm}.
\cref{fig:Usin} shows an accurate simulation, and contrasts it to with solution of~\cref{exForm}.
\end{enumerate}

\cref{sec:trust} details the accurate simulations of these four problems.
\cref{sec:dpat} develops a patch framework that accurately simulates all four of these archetype problems.
The framework extends the patch scheme, summarised by \cref{sec:pat}, by introducing a new so-called `double patch' that is placed over the shock.  

\section{Multiscale system simulation with the patch scheme} 
\label{sec:pat}

We present a basic implementation of the patch scheme in 1D space \cite[e.g.]{Gear03, Hyman2005, Samaey03b, Samaey04, Roberts06d}. 

\begin{figure}
    \centering
    \tikzsetfigurename{figs/patchsetup}%
\begin{tikzpicture}
\draw[->] (0,1) -- (11,1) node[right,below=0.3cm] {\textcolor{mtxt}{micro}scale} node[right,above=0.3cm] {\textcolor{Mtxt}{macro}scale} node[right] {$\textcolor{Mtxt}{X},\,\textcolor{mtxt}{x}$};
\filldraw[fill=mcol!20,draw=black] (30/22,1.1) rectangle (10-30/22,0.9);
\foreach \x in {3,...,19} {
	\draw[mtxt,thick] (10*\x/22,1.05)--(10*\x/22,0.95);
	}
\draw[mtxt,<->] (60/22,0.8)--(70/22,0.8) node[below,midway]{$d$};
\draw[Mcol,ultra thick] (30/22,1.2)--(30/22,0.8)  node[below] {\textcolor{mtxt}{$x_{j,-n}$}}; 
\filldraw[fill=Mcol!30, draw=Mcol, very thick] (5,1) circle (0.2cm) node[above,yshift=0.35cm] {\textcolor{Mtxt}{$X_j$}} node[below,yshift=-0.2cm] {\textcolor{mtxt}{$x_{j,0}$}}; 
\draw[mcol,<->] (5,1.35)--(190/22,1.35) node[above,midway]{\textcolor{mtxt}{$h$}};
\draw[Mcol, ultra thick] (10-30/22,1.2)--(10-30/22,0.8) node[below] {\textcolor{mtxt}{$x_{j,n}$}}; 
\end{tikzpicture}
    \caption{Schematic of the microscale detail of a single patch (\textcolor{mtxt}{coloured teal in the pdf version}).
Short vertical lines show the microscale mesh.
Three mesh-points on the patch connect to the macroscale (\textcolor{Mtxt}{orange in the pdf version}): the centre of the patch, and the two edge locations.
As shown in \cref{fig:pats} the simulation at each patch centre is coupled between patches to provide edge values to each patch.}
    \label{fig:patset}
\end{figure}
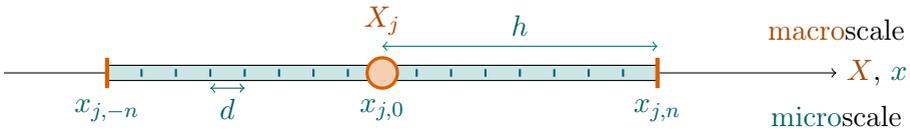

On the spatial domain, here \(-\pi<x<\pi\)\,, set equi-spaced macroscale nodes~$X_j$ on which the large-scale solution is to be computed, say $\Delta X_j = H$.
Centred on each of these macroscale nodes we place a `patch'. 
For all of our archetype problems, a patch consists of a finite difference discretisation of \pde~\cref{be} on a microscale mesh of $2n+1$~points with microscale equal-spacing~$\dx$ (\cref{fig:patset}).
Label the patches' microscale mesh points by $x_{j,i} = X_j+\dx i$ for microscale index $\mathcode`\,="213B i\in[-n,\ldots,-1,0,1,\ldots,n]$, and denote the patch half-width by $h:=n\dx$\,.
In general, these microscale parameters are to be set appropriate to the microscale.
In our archetype problems, they are set based on the diffusion coefficient~$\eps(u)$ in~\cref{be}.

\cref{fig:patset} illustrates these details. 
The microscale discretisation of the \pde~\eqref{be} is applied at the \(2n-1\)~interior points of each patch. 
In each patch, \cref{fig:patset} identifies three locations of interest on the macroscale: the centre `macroscale node' $x_{j,0}:= X_j$\,; and the two patch edges $x_{j,\pm n} = X_j \pm h$\,.
The remainder of this section describes the coupling of patches by the setting of the field edge values~$u_{j,\pm n}$ of every patch.
This coupling closes the patch simulation.

\begin{figure}
    \centering
    \tikzsetfigurename{figs/patches}%
\begin{tikzpicture}
\let\scriptstyle\relax
\tikzset{->-/.style={decoration={
  markings,
  mark=at position .45 with {\arrow{>}}},postaction={decorate}}}
  
\def\h{\textwidth/70};  
\def\hh{\h * 12}; 
\def\n{3}; 
\def\xx{\h * 30}; 
\draw[->](-\xx - 6*\h/2,1)--(\xx + 6*\h/2,1)node[right]{$X$};
\foreach \x in {-2,...,2} {\ifthenelse{\NOT \x=0}{
	\filldraw[mcol!20,draw=mtxt] (\x * \hh - \h * \n - \h/2,0.92)rectangle (\x * \hh + \h * \n + \h/2,1.08);
	\draw[Mcol,ultra thick] (\x * \hh - \h * \n - \h/2,1.15)--(\x * \hh - \h * \n - \h/2,0.85);
	\draw[Mcol,ultra thick] (\x * \hh + \h * \n + \h/2,1.15)--(\x * \hh + \h * \n + \h/2,0.85);
	\draw[excol,->-] (\hh*\x,1) .. controls ($(\hh*\x/2 - \h * \n - \h/4,1.5) + abs(\x)*(2*\h,1.5*\h)$) .. (- \h * \n - \h/2,1) ;
	\filldraw[fill=Mcol!30, draw=Mcol, thick] (\hh*\x,1) circle (0.9ex) node[below=1ex]{\textcolor{Mtxt}{$\scriptstyle X_{j\ifthenelse{0>\x}{\x}{+\x}}$}}
	 node[above=0.8ex]{{\textcolor{extxt}{$\scriptstyle U_{j\ifthenelse{0>\x}{\x}{+\x}}$}}}; 
	}{
	\filldraw[mcol!20,draw=mtxt] (\x * \hh - \h * \n - \h/2,0.92)rectangle (\x * \hh + \h * \n + \h/2,1.08);
	\draw[Mcol,very thick] (\x * \hh - \h * \n - \h/2,1.15)--(\x * \hh - \h * \n - \h/2,0.85);
	\draw[Mcol,very thick] (\x * \hh + \h * \n + \h/2,1.15)--(\x * \hh + \h * \n + \h/2,0.85);
	\draw[excol,->-] (\hh*\x,1) .. controls ( \hh*\x/2 - \h * \n/2 - \h/4,1.5) .. (- \h * \n - \h/2,1) ;
	 }
	 } 
	 \def\x{0} 
	 \filldraw[fill=Mcol!30, draw=Mcol, thick] (\hh*\x,1) circle (0.9ex) node[below=1ex]{\textcolor{Mtxt}{$\scriptstyle X_{j}$}}
	 node[above=2ex]{{\textcolor{extxt}{$\scriptstyle U_j$}}}; 
\draw[Mcol,<->,thick](-\hh, 0.1)--(0, 0.1) node [midway, below,yshift=-0.2] {\textcolor{Mtxt}{$H$}};			
\end{tikzpicture}
    \caption{Schematic showing inter-patch coupling to obtain the edge value at~$x=\textcolor{Mtxt}{X_j-h}$ (\cref{fig:patset} shows the internal structure of a patch).
Interpolating centre-patch field values~$\textcolor{extxt}{U_j}$ of the simulation at neighbouring macroscale nodes~$\textcolor{Mtxt}{X_j}$ provides all the edge values.
This figure illustrates the case of using \(\Gamma=2\) neighbours in each direction.
}
    \label{fig:pats}
\end{figure}
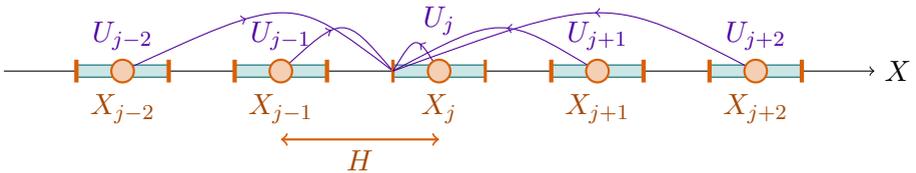

The key idea of the inter-patch coupling is to use the field values at macroscale nodes,~$\textcolor{extxt}{U_j}$, to set the field values~$u_{j,\pm n}$ at the patch edges~$x_{j,\pm n}$ (\cref{fig:pats}).
This is done by choosing a coupling order~\(\Gamma\), then for patch~$j$ use polynomial interpolation, to the patch edges~$x_{j,\pm n}$, of the centre-patch field values~$\textcolor{extxt}{U_j}$ in the~$2\Gamma+1$ nearest patches.
\cref{fig:pats} schematically illustrates this interpolation for one of the patch edges. 
In general, the patch edge values are the classic Lagrange interpolation
\begin{align} \label{patCoup}
u_{j,\pm n} := \sum_{i=j-\Gamma}^{j+\Gamma} \left(\prod_{k=j-\Gamma,\, k \ne i}^{j+\Gamma} \frac{ x_{j,\pm n} - X_k}{X_i - X_k}\right) U_i \,.
\end{align}
Other expressions, such as those using centred difference and mean operators \cite[]{Roberts06d, Roberts2011a, Cao2014a}, are equivalent in practice and may be used.

The patch scheme simulates a multiscale system on only the small fraction~$2h/H$ of the 1D domain (in~\(n\)D the fraction is~\((2h/H)^n\)).
If the characteristic microscale length scales are much smaller than the macroscale, $h\ll H$, then the patch scheme is much more efficient than computing over the full spatial domain.
Further, for a variety of problems, and as commented in the Introduction, the patch scheme has been proven to be generally consistent to the underlying system to order~$H^{2\Gamma}$ \cite[]{Roberts00a, Roberts2011a, Bunder2019c, Bunder2013b}. 

But such proven performance does not directly apply, and does not directly hold (see \cref{fig:stdPat} left), when a system has microscale dynamics that manifests as a macroscale shock.

\section{Shocks separate patch simulations into coupled macroscale systems}
\label{sec:dpat}

The success of the patch scheme relies on restricting, to the macroscale, almost all of the microscale information within each patch.
In the basic patch scheme (\cref{sec:pat}), we assume that the macroscale can be predicted appropriately from just \emph{one} characteristic of the dynamics in each patch, and so we may as well use the centre-patch value.
However, all four of the archetype problems~\cref{exStat,exForm,exUstat,exUform} violate this tenet in their shocks. 
A shock is a microscale transitional structure whose macroscale effects must be characterised by \emph{two} values: an average field value and the jump; or the two field values either side of the shock.
Thus a shock \emph{cannot} be resolved by patches as described by \cref{sec:dpat}---recall  that \cref{fig:stdPat}(left), in which patches cover an exorbitant one-third of the spatial domain, demonstrates the failure of the standard scheme to simulate~\cref{exForm} with its shock.
This section modifies and tests the patch scheme to overcome this problem.

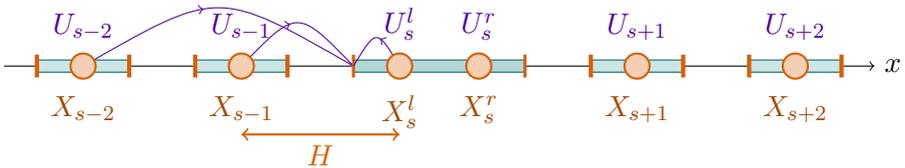
\begin{figure}
    \centering
    \tikzsetfigurename{figs/doublepatches}%
\begin{tikzpicture}
\let\scriptstyle\relax
\tikzset{->-/.style={decoration={
  markings,
  mark=at position .45 with {\arrow{>}}},postaction={decorate}}}
  
\def\h{\textwidth/70};  
\def\hh{\h * 12}; 
\def\n{3}; 
\def\xx{\h * 30}; 
\draw[->](-\xx - 3*\h,1)--(\xx + 3*\h,1)node[right]{$x$};
\foreach \x in {-2,-1,1,2} {
	\filldraw[mcol!20,draw=mtxt] ($  (\x * \hh - \h * \n - \h/2,0.92) + \x/abs(\x)*(\n*\h,0) $)rectangle ($ (\x * \hh + \h * \n + \h/2,1.08) + \x/abs(\x)*(\n*\h,0) $);
	\draw[Mcol,ultra thick] ($  (\x * \hh - \h * \n - \h/2,1.15) + \x/abs(\x)*(\n*\h,0) $)--($  (\x * \hh - \h * \n - \h/2,0.85) + \x/abs(\x)*(\n*\h,0) $);
	\draw[Mcol,ultra thick] ($ (\x * \hh + \h * \n + \h/2,1.15) + \x/abs(\x)*(\n*\h,0) $)--($ (\x * \hh + \h * \n + \h/2,0.85) + \x/abs(\x)*(\n*\h,0) $);
	\ifthenelse{0>\x}{
	\draw[excol,->-] ($ (\hh*\x,1) + \x/abs(\x)*(\n*\h,0) $) .. controls ($(\hh*\x/2 - \h * \n/2 - \h/4,1.5) + abs(\x)*(0,1.5*\h) + \x/abs(\x)*(\n*\h,0) $) .. (- 2*\h * \n - \h/2,1) ;
	}{}
	\filldraw[fill=Mcol!30, draw=Mcol, thick] ($ (\hh*\x,1) + \x/abs(\x)*(\n*\h,0) $) circle (1ex) node[below=1.5ex]{\textcolor{Mtxt}{$\scriptstyle X_{s\ifthenelse{0>\x}{\x}{+\x}}$}}
	 node[above=1.2ex]{{\textcolor{extxt}{$\scriptstyle U_{s\ifthenelse{0>\x}{\x}{+\x}}$}}}; 
	}
\def\x{0} 
	\filldraw[fill=mcol!30, draw=mtxt] (\x * \hh - 2*\h * \n - \h/2,0.92)rectangle (\x * \hh + 2*\h * \n + \h/2,1.08);
	\draw[Mcol,very thick] (\x * \hh - 2*\h * \n - \h/2,1.15)--(\x * \hh - 2*\h * \n - \h/2,0.85);
	\draw[Mcol,very thick] (\x * \hh + 2*\h * \n + \h/2,1.15)--(\x * \hh + 2*\h * \n + \h/2,0.85);
	\draw[excol,->-] (\hh*\x-\n*\h,1) .. controls ( \hh*\x/2 - 3*\h * \n/2 - \h/4,1.5) .. (- 2*\h * \n - \h/2,1) ;
	\filldraw[fill=Mcol!30, draw=Mcol, thick] (\hh*\x-\n*\h,1) circle (1ex) node[below=1.5ex]{\textcolor{Mtxt}{$\scriptstyle X^l_{s}$}}
	 node[above=1.2ex]{{\textcolor{extxt}{$\scriptstyle U^l_s$}}}; 
	 \filldraw[fill=Mcol!30, draw=Mcol, thick] (\hh*\x+\n*\h,1) circle (1ex) node[below=1.5ex]{\textcolor{Mtxt}{$\scriptstyle X^r_{s}$}}
	 node[above=1.2ex]{{\textcolor{extxt}{$\scriptstyle U^r_s$}}}; 
\draw[Mcol,<->,thick](-\hh-\n*\h, 0.1)--(0-\n*\h, 0.1) node [midway, below,yshift=-0.2] {$H$};			
\end{tikzpicture}
    \caption{This illustrates both a double patch (centre), with index \(j=s\), and the inter-patch coupling to obtain the edge value on the left side of the double patch. 
    The two `shock' nodes are labeled \textcolor{Mtxt}{$X^l_s$} and~\textcolor{Mtxt}{$X^r_s$}. 
    Inter-patch coupling is by the usual interpolation~\eqref{patCoup} (\cref{fig:pats}) except that the interpolation is adjusted so that it does not cross the double patch.}
    \label{fig:dpats}
\end{figure}

The key issue to be resolved is the communication of macroscale information across shocks---this communication must be done \emph{only} through the microscale simulation. 
We introduce a new object, called a `double patch' (\cref{fig:dpats} centre), which is to be placed over shocks.
The shock is resolved by the microscale discretisation inside the double patch, and the macroscale inter-patch coupling is altered as follows.
The double patch is so named, not because of its precise size, but because it contains \emph{two} macroscale nodes---called `shock nodes' for simplicity---that reflect the need for a shock to be characterised by \emph{two} macroscale variables.
If the double patch has patch index~\(j=s\)\,, then \(X_s^r\)~denotes the location of the right shock node, and \(X_s^l\)~denotes the left  (\cref{fig:dpats}).
The shock is to be located between the two shock nodes, and this constrains the minimum size of the double patch.  
Specifically, since the microscale is to resolve the full details of the shock, so the two `shock nodes' must lie on either side of the microscale transition layer (which appears on the macroscale as the shock).
Think of this as a simulation of two coupled macroscale domains, with the special double patch coupling the different macro-domains by providing the information to correctly evolve both the `left' and the `right' macro-domain. 

The macroscale inter-patch interpolation that couples patches is treated differently in the vicinity of the double patch. 
\begin{itemize}
\item Patch edges located to the left of the left shock node~\(X_s^l\) have their values determined by interpolation through the field values of \emph{only} those patch nodes in \(x\leq X_s^l\) (the `left' macro-domain, see \cref{fig:dpats}).
\item On the other side of the shock, patch edges located to the right of the right shock node~\(X_s^r\) have their values determined by interpolation through the field values of \emph{only} those patch nodes in \(x\geq X_s^r\) (the `right' macro-domain).
\end{itemize}
\cref{fig:dpats} illustrates these changes to the patch scheme. 
These changes are all that are needed to accurately simulate~\cref{exStat,exForm,exUstat,exUform}.

To be explicit, suppose there is a shock within patch~$s$.
Then there is no central macroscale node~$X_s$ in that patch.
Instead there are two nodes, $X_s^l$ and~$X_s^r$, to the left and right, respectively, of the shock's microscale transition layer.
Any patch~$j$ of sufficient distance from the shock, $|j-s| > \Gamma$, has patch edges specified by the usual Lagrange centred interpolation~\cref{patCoup}.
For patches close enough to the shock that they would normally couple to it, $|j-s| \le \Gamma$, then we implemented the option of a constant bandwidth adjusted interpolation.
To couple neighbours on the left of the shock, $s-\Gamma\leq j< s$, then the macroscale node of interest in the double patch is $X_s := X_s^l$, $U_s : = U_s^l$, and the edge values are
\newcommand\s{s}
\begin{align} \label{sCoupL}
u_{j,\pm n} := \sum_{i=j-\Gamma}^{\s} \left(\prod_{k=j-\Gamma,\, k \ne i}^{\s} \frac{ x_{j,\pm n} - X_k}{X_i - X_k}\right) U_i \,.
\end{align}
To couple neighbours on the right of the shock, $s< j\leq s+\Gamma$, set $X_s := X_s^r$, $U_s : = U_s^r$, and the edge values
\begin{align} \label{sCoupR}
u_{j,\pm n} := \sum_{i=\s}^{j+\Gamma} \left(\prod_{k=\s,\, k \ne i}^{j+\Gamma} \frac{ x_{j,\pm n} - X_k}{X_i - X_k}\right) U_i \,.
\end{align}
For the double patch itself, $j=s$\,, calculate~$u_{s,-n}$ with~\cref{sCoupL}, and~$u_{s,+n}$ with~\cref{sCoupR}.
Classical results \cite[e.g.]{Shoosmith75,Beyn79} guarantee that this constant bandwidth truncation, while it reduces the local order of consistency with the \pde, does not affect the global order of consistency.
Alternatively one may replace~\cref{sCoupL,sCoupR} with any good interpolation of the~$U_j$ that treats~$U_s^l$ and~$U_s^r$ as boundary values on the left and right sides of the shock, respectively, and also is of order~$2\Gamma$; for example, the standard asymmetric finite differences.

\cref{secces} discusses the motivation and theoretical support for the double patch. 
But first we simulate the archetype problems.

\subsection{Accurately simulate problems with the double patch}
\label{sPara}

\begin{figure}
    \centering
\begin{tabular}{cc}
    \includegraphics[]{M1P}
&   \includegraphics[]{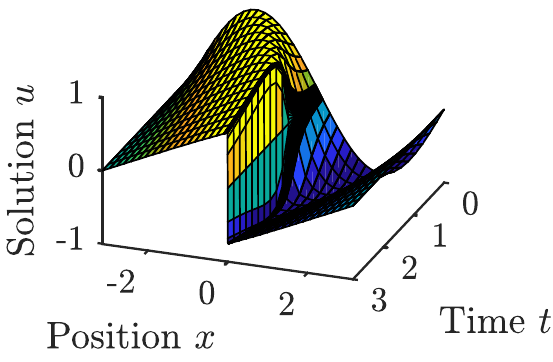}
\\    
    \includegraphics[]{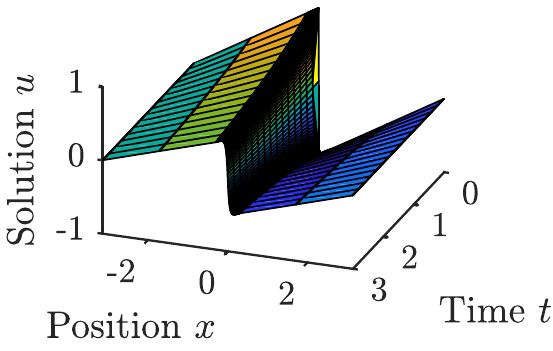}
&   \includegraphics[]{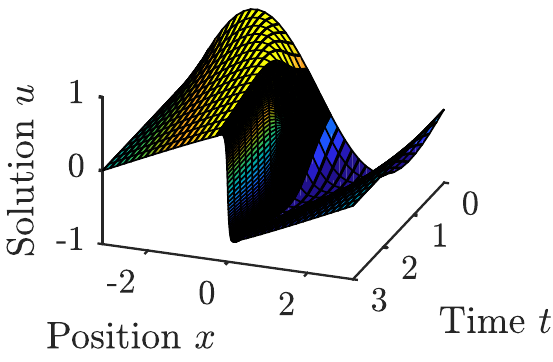}
\end{tabular}
    \caption{Numerical solution of~\cref{exStat,exForm,exUstat,exUform} by patches, using the double patch of \cref{sec:dpat} (compare to the accurate solutions in \cref{fig:tanh,fig:sin,fig:Utanh,fig:Usin}).
\cref{sPara} details each simulation.
The largest error, measured by comparing the simulation at macroscale nodes to the corresponding trusted solution, is~$0.0001$ for \cref{exStat} (top-left), $0.036$~for \cref{exForm} (top-right), $0.002$~for \cref{exUstat} (bottom-left), and $0.006$~for \cref{exUform} (bottom-right).  }
    \label{fig:P}
\end{figure}

All simulations reported here capture a macroscale shock inside a double patch centred at $x=0$, and spread patches uniformly over the remaining space.
\cref{fig:P} plots all simulations. 
The specific details are the following.%
\footnote{Although it is conceptually important that we use patches everywhere, to account for the practical situation in which there is not only a shock but also some multiscale character to the full problem (e.g., heterogeneous diffusion on the microscale), in these basic archetype problems there is no need to use patches except at the shock.
Essentially identical results are obtained by coupling the one double patch to a standard discretisation, with macroscale mesh spacing~\(H\), over the remaining spatial domain.
}
\begin{enumerate}
\item[\cref{exStat}] \emph{Simulate with few patches and a small double patch.} The double patch has width~\(0.05\) and contains \(25\)~mesh points. 
We discretise the remaining space into four patches, each of width~\(0.01\) and containing five mesh points.
We couple the patches with quadratic interpolation ($\Gamma=1$). 
This simulation displays features common with~\cref{exUstat}: the solution away from the shock is almost linear, and is accurately predicted with little effort. 
The double patch resolves the shock transition.

\item[\cref{exForm}] \emph{Simulate with many patches and a moderate double patch.} The double patch has width~\(0.2\) and contains \(100\)~mesh points.
We discretise the remaining space into \(34\)~patches, each of width~\(0.01\) and containing five mesh points, coupled with up to sixth order polynomial fits ($\Gamma=3$).
This simulation is both more expensive and less accurate than that for \cref{exStat} because of errors incurred \emph{during the formation} of the shock.
The largest error measured, at time~\(1.1\), is~$0.04$. 
However, outside the double patch the largest error over all times is~$0.006$.

\item[\cref{exUstat}] \emph{Simulate with few patches and a large double patch.} 
The double patch has width~\(0.6\) and contains \(180\)~mesh points. 
We discretise the remaining space into four patches, each of width~\(0.02\) and containing five mesh points, with quadratic coupling ($\Gamma=1$). 
Because of the nonlinearly enhanced diffusion, the field dissipates more broadly in the microscale, and so the double patch is wider to capture that behaviour---but the macroscale picture outside the double patch is unchanged. 

\item[\cref{exUform}] \emph{Simulate with many patches and a large double patch.} 
The double patch has width~\(0.6\) and contains \(180\)~mesh points.
We discretise the remaining space into \(34\)~patches, each of width~\(0.02\) and containing four mesh points, coupled with up to sixth order polynomials ($\Gamma=3$). 
This simulation is the largest of all archetype problems, covering~20\% of the spatial domain. 
From the sinusoidal initial condition of~\cref{exForm} it also displays a higher error in the initial formation of the shock, and from the variable diffusion coefficient~$\eps(u)$ of~\cref{exUstat} it inherits the large double patch that is required to accurately simulate the diffusing shock edges. 
\end{enumerate}

\subsection{Consistency of the extended scheme} 
\label{secces}

The patch scheme aims to use a microscale simulator that a user certifies is accurate enough for the purposes at hand.
The scheme couples together the microscale simulations in a way that is consistent to order~$H^{2\Gamma}$, where $H$~is the spacing between patches (\cref{fig:pats}).
We carefully designed a modification to the patch scheme that allows a shock to be hidden within a particular patch.
Notionally then, the overall error resembles that of any patch scheme,~$H^{2\Gamma}$, \emph{plus} any considerations from the placement of the shock in the double patch and the coupling near the shock.%
\footnote{We assume that a user's microscale simulation resolves the physical details and dynamics within the shock to an adequate accuracy.
}
We now establish one way to control these additional sources of error. 

That the coupling near the shock is appropriate was established below~\cref{sCoupR}: our proposed coupling is consistent on the macroscale using results by \cite{Shoosmith75,Beyn79}.
The remaining errors are incurred by the location of the shock and of the two shock nodes. 

Let us briefly consider~\cref{exStat} for $x\ge0$ to the right of the shock: after an initial transient the solution is $u(x,t) \approx (1-e^{-x/\eps} ) A(t) (x-\pi)$.
That is, the solution consists of an exponentially decaying term resolving the inner shock transition, and a linear outer solution $u = A(t)(x-\pi)$ away from the shock. 
Simulating with the novel patch scheme with a shock patch of index $j=s$ centred at $x=0$, consistency to the same order as a typical patch scheme is assured so long as $\exp\big(-X_s^r/\eps\big) \ll H^{2\Gamma}$, and by symmetry a similar condition holds for~$X_s^l$. 
That is, the general rule that we require is that the shock nodes be sufficiently widely spaced that the microscale structure of the shock influence the field values~$U_s^l$ and~$U_s^r$ by no more than the desired accuracy of consistency. 
For example, for three significant digit accuracy, and since \(e^{-7}\approx 0.001\) we would require that \(X_s^r\approx (7-2\Gamma\log H)\eps\).
That is, typically we expect to need the double patch's nodes to be separated by some relatively moderate multiple of~\(\eps\).

\section{Discussion}
Here we extended the patch scheme to cater for the scenario in which a multiscale system contains one or more shocks, where a `shock' is defined as a localised microscale feature that transitions between two macroscale domains with quantitatively different solutions (e.g., a crack, dislocation, grain boundary).
The innovation is to resolve full microscale details of each shock inside a so-called `double patch'.
Then near each of the two edges of the double patch, we define left and right nodes: these nodes provide field values for the macroscale predictions on the corresponding sides of the double patch, and which then  contribute to coupling the patches together. 
This careful treatment of the shocks allows accurate, rapid simulations to be performed in the usual multiscale picture of the patch scheme.
In keeping with the philosophy of the patch scheme, we do not require any detailed information on the shocks, or whether or not some quantities are conserved or not, only that the shock location is approximately known.
 
In ongoing research, we are extending the double patch concept to moving patches, and allow merging patches, by adapting techniques developed for moving meshes \cite[e.g.]{Budd2009}.
Such movement and merging will empower patches to adaptively track emerging shocks, and significantly reduce the computational burden of simulating problems like~\cref{exForm,exUform}.
However, for problems with a shock existing in the initial state, like~\cref{exStat,exUstat}, the methods in this article will often be sufficient.

Future research will look at the issues involved in the more difficult task of predicting shocks in multi-D by extending the multi-D patch scheme \cite[e.g.]{Roberts2011a, Bunder2019c}.

\paragraph{Acknowledgement}
This research was funded by the Australian Research Council under grants DP150102385 and DP180100050.  
The work of I.G.K. was also partially supported by the DARPA PAI program.

\appendix

\section{Trusted solutions to archetype problems}
\label{sec:trust}
In order to study the effectiveness of the proposed approach, we computed `trusted solutions', close approximations of the exact solution, to each archetype problem. 
These are described in this appendix.

\subsection{Simulate \cref{exStat,exForm} with a numerical quadrature of the exact solution}

The well-known exact solution to Burgers' \pde~\cref{be} uses the Cole--Hopf transformation \cite[e.g.]{Whitham74} to relate solutions to those of the linear diffusion \pde.
The diffusion \pde\ is solved exactly and then, reversing the transformation, the exact solution for Burgers' \pde\ at every position and time is
\begin{align}
    \label{beEx}
    u(x,t) &= \frac{\int_{-\infty}^\infty (x-y) \exp\left[-\frac{(x-y)^2}{4\eps t} - \frac{1}{2\eps}\int_0^y u_0(z)\,dz\right]dy}
    {t\int_{-\infty}^\infty \exp\left[-\frac{(x-y)^2}{4\eps t} - \frac{1}{2\eps}\int_0^y u_0(z)\,dz\right]dy} \,.
\end{align}
We also exploited that solutions to both~\cref{exStat,exForm} have $u=0$ at $x=\pm \pi$ for all time, since the numerical quadrature~\cref{beEx} does not straightforwardly cater for boundary conditions. 

Notionally, we could compute the numerator and denominator of~\cref{beEx} using numerical quadratures.
But in practice the integrals are difficult to compute at small values of~$\eps$ because of rapid variations in the integrands.
A direct evaluation via quadratures of~\cref{beEx} using \verb|integral()| in \Matlab\ is unstable around $x=0$ with even the moderate choice of diffusion $\eps=0.01$.
Instead we modified~\cref{beEx}, summarised in~\cref{beTr}, to allow us to compute a trusted solution for all $\eps \ge 10^{-4}$.
The integrals~\cref{beEx} are adjusted in two ways in order to compute them.
First, at each~$x$ we identify the value of~$y$ that maximises the argument of the exponentials in~\cref{beEx}, $v(x,y) = -\frac{(x-y)^2}{4 t} - \frac{1}{2}\int_0^y u_0(z)\,dz$; that is, we identify $y^*(x) = \argmax_y v(x,y)$.
We change the limits of the integral to only compute from $y^*(x)-\tol$ to $y^*(x)+\tol$, for some tolerance $\tol>0$.
The value of each integrand decreases rapidly around~$y^*(x)$, scaling roughly like $\exp[-(y-y^*(x))^2/\eps]$, so this approximation is extremely accurate at even small values $\tol=100\epsilon$; we chose $\tol=5$ for safety.
The second adjustment to~\cref{beEx} is that we scale both numerator and denominator by $C(x) = v(x,y^*(x))$.
This does not algebraically affect the solution, and is done to avoid round-off errors.
Summarising, numerical quadrature determines the trusted solution via
\begin{align}
    \label{beTr}
    u(x,t) &= \frac{\int_{y^*(x)-5}^{y^*(x)+5} (x-y) \exp\left[\frac1{\eps}\left\{-\frac{(x-y)^2}{4 t} - \frac{1}{2}\int_0^y u_0(z)\,dz - C(x)\right\}\right]dy}
    {t\int_{y^*(x)-5}^{y^*(x)+5} \exp\left[\frac1{\eps}\left\{-\frac{(x-y)^2}{4 t} - \frac{1}{2}\int_0^y u_0(z)\,dz - C(x)\right\}\right]dy} \,.
\end{align}

\subsection{Brute force approaches simulate \cref{exUstat,exUform}}
These problems cannot be solved with numerical quadratures because the exact algebraic solution is not available.
Instead we discretise the modified Burgers' \pde~\cref{be} by finite differences,
\[
\frac{du_i}{dt} =  \frac{1}{d^2}(\eps_1+\eps_2|u_i|)(u_{i+1}-2u_i+u_{i-1}) - \frac{u_i}{2d}(u_{i+1} - u_{i-1})
\]
on a fine grid of \(1600\)~points between~$-\pi$ and~$\pi$, with spacing $d=0.00375$ between each grid point.
This discretisation is simulated with a fast-time step of~$d^2/2$ up until the desired final time---this consumes vastly more computational time than the patch scheme and is done only to determine errors.
Lastly, the fine solution is interpolated to desired spatial locations and times.

%


@article{Kevrekidis03b,
  author =        {Kevrekidis, I.~G. and Gear, C.~W. and Hyman, J.~M. and
                   Kevrekidis, P.~G. and Runborg, O. and
                   Theodoropoulos, K.},
  journal =       {Comm. Math. Sciences},
  pages =         {715--762},
  title =         {Equation-free, coarse-grained multiscale computation:
                   enabling microscopic simulators to perform system
                   level tasks},
  volume =        {1},
  year =          {2003},
}

@article{Kevrekidis09a,
  author =        {Kevrekidis, Ioannis~G. and Samaey, Giovanni},
  journal =       {Annu. Rev. Phys. Chem.},
  pages =         {321--44},
  title =         {Equation-Free Multiscale Computation: Algorithms and
                   Applications},
  volume =        {60},
  year =          {2009},
  doi =           {10.1146/annurev.physchem.59.032607.093610},
}

@article{Sieber2018,
  author =        {Sieber, J. and Marschler, C. and Starke, J.},
  journal =       {SIAM Journal on Applied Dynamical Systems},
  month =         jan,
  number =        {4},
  pages =         {2574--2614},
  title =         {Convergence of {Equation}-{Free} {Methods} in the
                   {Case} of {Finite} {Time} {Scale} {Separation} with
                   {Application} to {Deterministic} and {Stochastic}
                   {Systems}},
  volume =        {17},
  year =          {2018},
  doi =           {10.1137/17M1126084},
}

@article{Gear02b,
  author =        {Gear, C. W. and Kevrekidis, Ioannis G.},
  journal =       {SIAM Journal on Scientific Computing},
  number =        {4},
  pages =         {1091--1106},
  publisher =     {SIAM},
  title =         {Projective Methods for Stiff Differential Equations:
                   Problems with Gaps in Their Eigenvalue Spectrum},
  volume =        {24},
  year =          {2003},
  doi =           {10.1137/S1064827501388157},
  url =           {http://link.aip.org/link/?SCE/24/1091/1},
}

@article{Erban2006,
  author =        {Erban, Radek and Kevrekidis, Ioannis G. and
                   Othmer, Hans G.},
  journal =       {Physica D: Nonlinear Phenomena},
  number =        {1},
  pages =         {1--24},
  title =         {An equation-free computational approach for
                   extracting population-level behavior from
                   individual-based models of biological dispersal},
  volume =        {215},
  year =          {2006},
  doi =           {10.1016/j.physd.2006.01.008},
  issn =          {0167-2789},
  url =           {http://www.sciencedirect.com/science/article/B6TVK-4JDVNSP-
                  1/2/f31e03e0a32cfcb2a811f41ed6a8dfc6},
}

@article{Givon06,
  author =        {Givon, Dror and Kevrekidis, Ioannis~G. and
                   Kupferman, Raz},
  journal =       {Comm. Math. Sci.},
  number =        {4},
  pages =         {707--729},
  title =         {Strong convergence of projective integration schemes
                   for singularly perturbed stochastic differential
                   systems},
  volume =        {4},
  year =          {2006},
}

@article{Gear03,
  author =        {Gear, C.~W. and Li, Ju and Kevrekidis, I.~G.},
  journal =       {Phys. Lett.~A},
  pages =         {190--195},
  title =         {The gap-tooth method in particle simulations},
  volume =        {316},
  year =          {2003},
  doi =           {10.1016/j.physleta.2003.07.004},
}

@article{Hyman2005,
  author =        {Hyman, James M.},
  journal =       {Computing in Science \& Engineering},
  number =        {3},
  pages =         {47--53},
  title =         {Patch Dynamics for Multiscale Problems},
  volume =        {7},
  year =          {2005},
  doi =           {10.1109/MCSE.2005.57},
  url =           {http://scitation.aip.org/content/aip/journal/cise/7/3/
                  10.1109/MCSE.2005.57},
}

@article{Samaey03b,
  author =        {Samaey, G. and Kevrekidis, I.~G. and Roose, D.},
  journal =       {Multiscale Modeling and Simulation},
  pages =         {278--306},
  title =         {The gap-tooth scheme for homogenization problems},
  volume =        {4},
  year =          {2005},
  doi =           {10.1137/030602046},
}

@article{Samaey04,
  author =        {Samaey, Giovanni and Roose, Dirk and
                   Kevrekidis, Ioannis~G.},
  journal =       {J.~Comput Phys.},
  pages =         {264--287},
  title =         {Patch dynamics with buffers for homogenization
                   problems},
  volume =        {213},
  year =          {2006},
  doi =           {10.1016/j.jcp.2005.08.010},
}

@article{Roberts06d,
  author =        {Roberts, A.~J. and Kevrekidis, I.~G.},
  journal =       {SIAM J.~Scientific Computing},
  number =        {4},
  pages =         {1495--1510},
  title =         {General tooth boundary conditions for equation free
                   modelling},
  volume =        {29},
  year =          {2007},
  doi =           {10.1137/060654554},
}

@article{Takayama04,
  author =        {Takayama, Kazuyoshi and Saito, Tsutomu},
  journal =       {Annu. Rev. Fluid Mechanics},
  pages =         {347--379},
  title =         {Shock wave/geophysical and medical applications},
  volume =        {36},
  year =          {2004},
  doi =           {10.1146/annurev.fluid.36.050802.121954},
}

@book{Guenter2011,
  author =        {G\"uenter Hofstetter and G\"uenther Meschke},
  note =          {\doi{10.1007/978-3-7091-0897-0}},
  publisher =     {Springer},
  series =        {{CISM} International Centre for Mechanical Sciences},
  title =         {Numerical Modeling of Concrete Cracking},
  year =          {2011},
}

@article{Luskin2013,
  author =        {Luskin, Mitchell and Ortner, Christoph},
  journal =       {Acta Numerica},
  month =         {5},
  pages =         {397--508},
  title =         {Atomistic-to-continuum coupling},
  volume =        {22},
  year =          {2013},
  doi =           {10.1017/S0962492913000068},
  issn =          {1474-0508},
  url =           {http://journals.cambridge.org/article_S0962492913000068},
}

@article{Tadmor1996,
  author =        {Tadmor, E.~B. and Ortiz, M. and Phillips, R.},
  journal =       {Philosophical Magazine~A},
  number =        {6},
  pages =         {1529--1563},
  title =         {Quasicontinuum analysis of defects in solids},
  volume =        {73},
  year =          {1996},
  doi =           {10.1080/01418619608243000},
}

@incollection{Samaey08,
  author =        {Samaey, G. and Roberts, A.~J. and Kevrekidis, I.~G.},
  booktitle =     {Multiscale methods: bridging the scales in science
                   and engineering},
  chapter =       {8},
  editor =        {Fish, Jacob},
  pages =         {216--246},
  publisher =     {Oxford University Press},
  title =         {Equation-free computation: an overview of patch
                   dynamics},
  year =          {2010},
  isbn =          {978-0-19-923385-4},
}

@book{Roberts2014a,
  author =        {Roberts, A.~J.},
  month =         {jan},
  publisher =     {SIAM, Philadelphia},
  title =         {Model emergent dynamics in complex systems},
  year =          {2015},
  isbn =          {9781611973556},
  url =           {http://bookstore.siam.org/mm20/},
}

@article{Bunder2013b,
  author =        {Bunder, J.~E. and Roberts, A.~J. and
                   Kevrekidis, I.~G.},
  journal =       {J.~Computational Physics},
  month =         {may},
  pages =         {154-174},
  title =         {Good coupling for the multiscale patch scheme on
                   systems with microscale heterogeneity},
  volume =        {337},
  year =          {2017},
  doi =           {10.1016/j.jcp.2017.02.004},
}

@article{Geers2010,
  author =        {Geers, M. G. D. and Kouznetsova, V. G. and
                   Brekelmans, W. A. M.},
  journal =       {Journal of Computational and Applied Mathematics},
  month =         aug,
  number =        {7},
  pages =         {2175--2182},
  series =        {Fourth {International} {Conference} on {Advanced}
                   {Computational} {Methods} in {Engineering} ({ACOMEN}
                   2008)},
  title =         {Multi-scale computational homogenization: {Trends}
                   and challenges},
  volume =        {234},
  year =          {2010},
  doi =           {10.1016/j.cam.2009.08.077},
  issn =          {0377-0427},
}

@article{Saeb2016,
  author =        {Saeb, Saba and Steinmann, Paul and Javili, Ali},
  journal =       {Applied Mechanics Reviews},
  month =         sep,
  number =        {5},
  title =         {Aspects of {Computational} {Homogenization} at
                   {Finite} {Deformations}: {A} {Unifying} {Review}
                   {From} {Reuss}' to {Voigt}'s {Bound}},
  volume =        {68},
  year =          {2016},
  doi =           {10.1115/1.4034024},
  issn =          {0003-6900},
}

@incollection{Geers2017,
  author =        {Geers, Marc G. D. and Kouznetsova, Varvara G. and
                   Matou{\v s}, Karel and Yvonnet, Julien},
  booktitle =     {Encyclopedia of Computational Mechanics, Second
                   Edition},
  pages =         {1--34},
  publisher =     {Wiley},
  title =         {Homogenization Methods and Multiscale Modeling:
                   Nonlinear Problems},
  year =          {2017},
  doi =           {10.1002/9781119176817.ecm2107},
  isbn =          {978-1-119-17681-7},
  url =           {https://onlinelibrary.wiley.com/doi/abs/10.1002/
                  9781119176817.ecm2107},
}

@incollection{Craster2015,
  author =        {Craster, Richard V.},
  booktitle =     {Analytic Methods in Interdisciplinary Applications},
  editor =        {Mityushev, V.~V. and Ruzhansky, M.},
  pages =         {41--50},
  publisher =     {Springer},
  series =        {Springer Proceedings in Mathematics and Statistics},
  title =         {Dynamic Homogenization},
  volume =        {116},
  year =          {2015},
  doi =           {10.1007/978-3-319-12148-2_3},
}

@article{Owhadi2015,
  author =        {Owhadi, Houman},
  journal =       {Multiscale Modeling \& Simulation},
  number =        {3},
  pages =         {812--828},
  title =         {Bayesian Numerical Homogenization},
  volume =        {13},
  year =          {2015},
  doi =           {10.1137/140974596},
  url =           {http://dx.doi.org/10.1137/140974596},
}

@techreport{Peterseim2019,
  author =        {Daniel Peterseim},
  institution =   {AMSI Winter School on Computational Modeling of
                   Heterogeneous Media,
  \url{https://ws.amsi.org.au/wp-content/uploads/sites/70/2019/06/numhomamsi2019.pdf}},
  month =         {6},
  title =         {Numerical Homogenization beyond Scale Separation and
                   Periodicity},
  year =          {2019},
}

@article{Maier2019,
  author =        {Maier, Roland and Peterseim, Daniel},
  journal =       {BIT Numerical Mathematics},
  month =         jun,
  number =        {2},
  pages =         {443--462},
  title =         {Explicit computational wave propagation in
                   micro-heterogeneous media},
  volume =        {59},
  year =          {2019},
  doi =           {10.1007/s10543-018-0735-8},
  issn =          {1572-9125},
}

@article{Roberts2011a,
  author =        {Roberts, A.~J. and MacKenzie, Tony and
                   Bunder, Judith},
  journal =       {J.~Engineering Mathematics},
  number =        {1},
  pages =         {175--207},
  title =         {A dynamical systems approach to simulating macroscale
                   spatial dynamics in multiple dimensions},
  volume =        {86},
  year =          {2014},
  doi =           {10.1007/s10665-013-9653-6},
  url =           {http://arxiv.org/abs/1103.1187},
}

@article{Cao2014a,
  author =        {Cao, Meng and Roberts, A. J.},
  journal =       {IMA J.~Applied Maths.},
  number =        {2},
  pages =         {228--254},
  title =         {Multiscale modelling couples patches of nonlinear
                   wave-like simulations},
  volume =        {81},
  year =          {2016},
  doi =           {10.1093/imamat/hxv034},
}

@book{Huang10,
  author =        {Huang, Weizhang and Russell, Robert D},
  publisher =     {Springer Science \& Business Media},
  title =         {Adaptive moving mesh methods},
  volume =        {174},
  year =          {2010},
}

@article{JakemanArchibaldXiu11,
  author =        {Jakeman, John D and Archibald, Richard and
                   Xiu, Dongbin},
  journal =       {Journal of Computational Physics},
  number =        {10},
  pages =         {3977--3997},
  publisher =     {Elsevier},
  title =         {Characterization of discontinuities in
                   high-dimensional stochastic problems on adaptive
                   sparse grids},
  volume =        {230},
  year =          {2011},
}

@article{Alotaibi2017a,
  author =        {Alotaibi, Hammad and Cox, Barry and Roberts, A. J.},
  institution =   {\url{http://arxiv.org/abs/1703.00204}},
  journal =       {ANZIAM~J.},
  number =        {3},
  pages =         {313--334},
  title =         {Couple microscale periodic patches to simulate
                   macroscale emergent dynamics},
  volume =        {59},
  year =          {2018},
  doi =           {10.1017/S1446181117000396},
}

@article{Mones2014,
  author =        {Mones, Enys and Araujo, Nuno A. M. and Vicsek, Tamas and
                   Herrmann, Hans J.},
  journal =       {Sci. Rep.},
  pages =         {4949},
  title =         {Shock waves on complex network},
  volume =        {4},
  year =          {2014},
  doi =           {10.1038/srep04949},
}

@article{Roberts00a,
  author =        {Roberts, A.~J.},
  journal =       {Mathematics of Computation},
  pages =         {247--262},
  title =         {A holistic finite difference approach models linear
                   dynamics consistently},
  volume =        {72},
  year =          {2003},
  doi =           {10.1090/S0025-5718-02-01448-5},
  url =           {http://www.ams.org/mcom/2003-72-241/S0025-5718-02-01448-5},
}

@techreport{Bunder2019c,
  author =        {J.E. Bunder and J. Divahar and Ioannis G. Kevrekidis and
                   Trent W. Mattner and A.J. Roberts},
  institution =   {\url{https://arxiv.org/abs/1912.07815}},
  month =         {dec},
  title =         {Large-scale simulation of shallow water waves with
                   computation only on small staggered patches},
  year =          {2019},
}

@incollection{Shoosmith75,
  author =        {Shoosmith, John N.},
  booktitle =     {Numerical {Solutions} of {Boundary} {Value}
                   {Problems} for {Ordinary} {Differential} {Equations}},
  editor =        {Aziz, A. K.},
  month =         jan,
  pages =         {355--369},
  publisher =     {Academic Press},
  title =         {A high-order finite-difference method for the
                   solution of two-point boundary-value problems on a
                   uniform mesh},
  year =          {1975},
  doi =           {10.1016/B978-0-12-068660-5.50020-X},
  isbn =          {978-0-12-068660-5},
}

@article{Beyn79,
  author =        {Beyn, Wolf-J{\"u}rgen},
  journal =       {Mathematics of Computation},
  number =        {148},
  pages =         {1213--1228},
  title =         {The exact order of convergence for finite difference
                   approximations to ordinary boundary value problems},
  volume =        {33},
  year =          {1979},
  doi =           {10.1090/S0025-5718-1979-0537966-2},
  issn =          {0025-5718, 1088-6842},
}

@article{Budd2009,
  author =        {Budd, Chris~J. and Huang, Weizhang and
                   Russell, Robert~D.},
  journal =       {Acta Numerica},
  pages =         {111--241},
  title =         {Adaptivity with moving grids},
  volume =        {18},
  year =          {2009},
  doi =           {10.1017/S0962492906400015},
}

@book{Whitham74,
  author =        {G.B. Whitham},
  publisher =     {John Wiley \& Sons},
  title =         {Linear and nonlinear waves},
  year =          {1974},
}
\end{document}